# Dynamic of the extended Mandelbrot's equation

Marek Berezowski

Silesian University of Technology, Institute of Mathematics
Gliwice, Poland
E-mail: marek.berezowski@polsl.pl

## Abstract

The paper deals with the mathematical–numerical analysis of the Mandelbrot equation extended by the dynamic continuous term. The possibilities of generation of fractal patterns with the mathematical form, defined in such a manner, were investigated. The decay of fractal structure with time was demonstrated.

## 1. Introduction

There are many works devoted to the Mandelbrot's equation and set. In order to recall it, Mandelbrot's equation constitutes a recurrent form

$$z_{k+1} = z_k^2 + c \qquad (1)$$

in which the variable $z$ and constant $c$ are complex numbers. In the result of recurrence Eq. (1) one obtains, for $z_0 = 0$, the sequence $z_k$, which is convergent or divergent, depending on the value of constant $c$. Confirming the convergence of this sequence by means of the point on the complex plane $c(c_r, c_i)$, results in the fractal structure, named—after its discoverer—Mandelbrot's fractal (Fig. 1). In spite of the fact that this structure is extremely complicated, it is possible to determine analytically its principal fragments. Viz., it may easily be shown that the edge of its largest area is a set of Hopf bifurcation points satisfying the condition $|z| = 1/2$. This set may be described by the parametric equation of the form

$$c_H = \frac{1}{2}e^{i\varphi} - \frac{1}{4}e^{2i\varphi} \qquad (2)$$

where parameter $\phi$ is an argument of the complex variable $z$ (full line in Fig. 1). Also quite easily one can determine the edge of solutions of two-period oscillations

$$z_{k+2} = (z_k^2 + c)^2 + c \qquad (3)$$

which satisfies the condition $|y| = 1/4$, where

$$y = z(z^2 + c) \qquad (4)$$

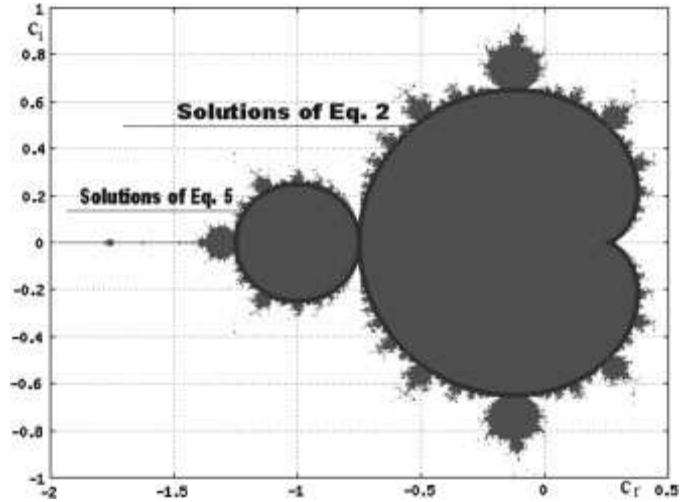
Fig. 1. Mandelbrot's fractal.

The parametric equation of this equation has the form

$$c_{2H} = \frac{1}{4}e^{i\varphi} - 1 \qquad (5)$$

where parameter $\phi$ is the argument of complex variable $y$. In a similar way one may also determine the edge of three-period solutions, for which the corresponding complex variable fulfils the condition $|x| = 1/8$, where

$$x = z(z^2+c)[(z^2+c)^2+c]. \qquad (6)$$

The relevant discussion was carried out a.o. in [1] and [2].

## 2. Transition from a discrete equation to a continuous equation

Let us assume that the general recurrent form

$$z_{k+1} = F(z_k) \qquad (7)$$

constitutes a mathematical model of some dynamic system with delay. Hence, such a system may be described by an equivalent continuous equation of the form

$$z(t) = F[z(t-t_0)] \qquad (8)$$

where $t_0$ is a certain delay time. If, additionally, given dynamic system is characterized by a significant inertion, the corresponding derivative with respect to time must appear in Eq. (8) [3] and [4], which leads to the relationship of the form

$$\sigma \frac{dz(t)}{dt} + z(t) = F[z(t-t_0)] \qquad (9)$$

Let us assume in our considerations that the mathematical relationship of Mandelbrot is connected with a certain continuous dynamic system comprising delay and inertion, the letter characterized by the dynamic capacity $\sigma$. The Mandelbrot's equation (1) extended in such a way has the form:

$$\sigma \frac{dz(t)}{dt} + z(t) = z^2(t - t_0) + c \tag{10}$$

Introducing the normalized time $\tau = \frac{t}{\sigma}$ one obtains the final equation

$$\frac{dz(\tau)}{d\tau} + z(\tau) = z^2(\tau - \tau_0) + c \tag{11}$$

In order to analyze the behavior of solutions of this equation, use has been made of the above-mentioned Mandelbrot's algorithm, which generates the fractal set. To this aim there was constructed the map of stable and unstable solution of Eq. (11) as a function of successive values of complex constant $c$ under the assumption $z(0) = 0$. For $\tau_0 = 0$ Eq. (11) reduces to the differential form without delay. In the consequence its solution cannot be of chaotic character and therefore the corresponding region cannot display the fractal character. Another extreme case is $\tau_0 = \infty$, which means that the form Eq. (11) is equivalent with the recurrent form Eq. (1) of Mandelbrot. In the work presented some exemplary computations were carried out under the assumption $\tau_0 = 10$; the corresponding image is shown in Fig. 2. It may be seen that this image is very similar to the original Mandelbrot's fractal form Fig. 1, which seems to be a quite natural result. This topic will be discussed later.

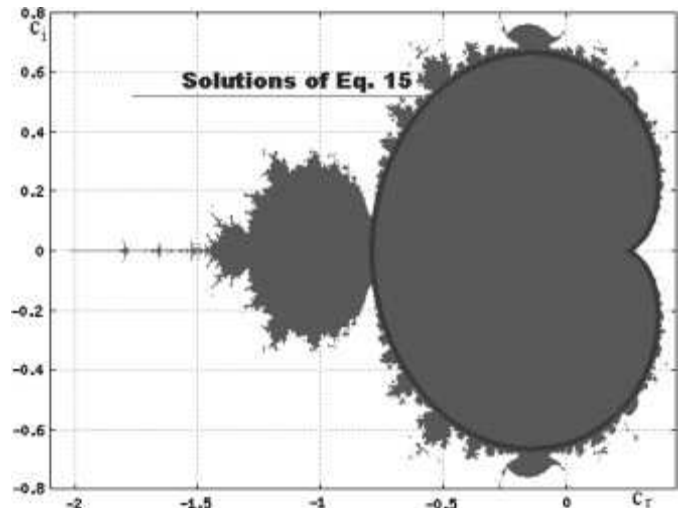

Fig. 2. Image for small time simulation.

Similarly as in the case of the logistic model Eq. (1), taking into consideration Eq. (11), one can obtain here the analytical relationship, which determines the edge of stable stationary solutions of the model (full line in Fig. 2). Namely, the linear approximation of Eq. (11), determined in the vicinity of stationary solution, is

$$\frac{du}{d\tau} + u = 2z_s u(\tau - \tau_0) \tag{12}$$

where $u = z - z_s$, $z_s$ being the stationary solution of relation Eq. (11). The characteristic equation of relationship Eq. (12), written in the frequency form, is

$$1 - \frac{2|z_s|e^{i\varphi}}{\sqrt{1+\omega^2}} e^{-i\operatorname{arctg}\omega} e^{-i\omega\tau_0} = 0 \tag{13}$$

Using Nyquist criterion one obtains the boundary of stability of Eq. (11) in the form

$$|z_s| = \frac{\sqrt{1+\omega^2}}{2} \tag{14}$$

Making use of the definition of the stationary state $(dz/d\tau = 0)$ one obtains finally the parametric relationship determining the edge of stationary solutions of Eq. (11) in the form

$$c_H = \frac{\sqrt{1+\omega^2}}{2} e^{i\varphi} - \frac{1+\omega^2}{4} e^{2i\varphi} \tag{15}$$

where parameter $\phi$ is an argument of complex variable $z_s$; the circular frequency $\omega$ is determined from the relation

$$\omega\tau_0 + \operatorname{arctg}\omega = \varphi \tag{16}$$

The regions outside the main part of Mandelbrot's set refer to the multiperiodic oscillations [1] and [2]. It comes out that during the formation of the structure as in Fig. 2 the magnitude of these external regions diminishes with the increase of the time of simulation of Eq. (11) (Fig. 3). This means that for the infinitely long simulation time all regions vanish, except the main one. This results from the fact that the oscillation solutions of Eq. (11) (except the case $c_i = 0$) are unstable. So, the geometric structure generated according to the above-mentioned algorithm consists exclusively of the main region, which characterizes the stable stationary solutions of this equation and of a single line, which starts to the left of this region from the point of the coordinates $(-0.791;0)$ (Fig. 4). This line, as the only exception, constitutes the set of stable periodic and chaotic solutions of Eq. (11), as shown in the Feigenbaum's diagram of steady states (Fig. 5). Similar results are obtained for any arbitrary delay time $\tau_0 < \infty$.

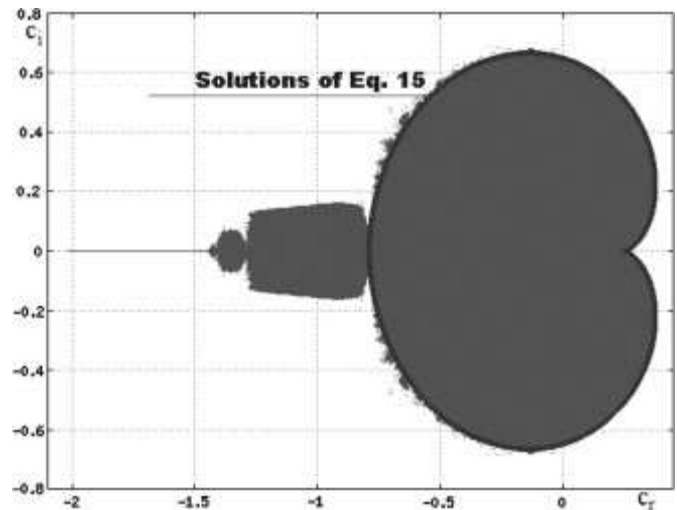
Fig. 3. Image for large time simulation.

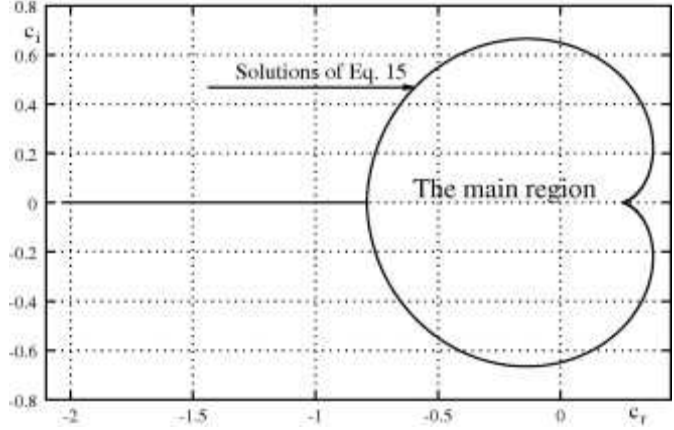

Fig. 4. Image for infinite time simulation.

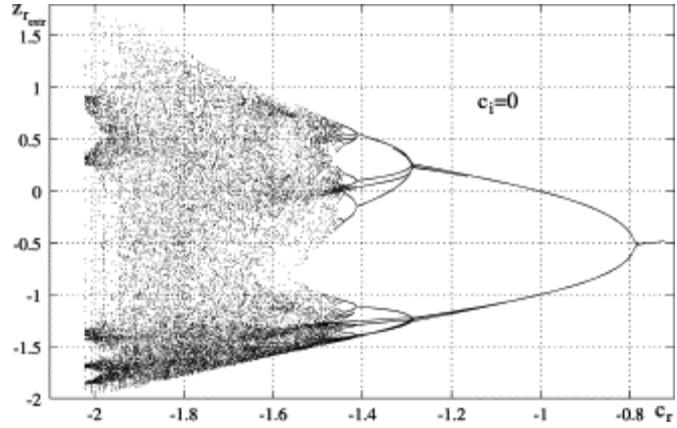

Fig. 5. Feigenbaum diagram.

So, the Mandelbrot's equation (1), as a particular case of the form Eq. (11) is an exception, which generates the fractal structure. Viz., if in this equation the inertial term appears, the fractal vanishes. This phenomenon may be investigated, analyzing the time series. First, however, let us discuss more precisely the edge of the region from Fig. 5, i.e. the set of points of Hopf bifurcation. To each of these points a certain value of the angle $\phi$ corresponds, and thus—the circular frequency $\omega$. The phase displacement, due to the right-hand side of Eq. (11) depends solely on $\phi$ and $\tau_0$. It results from Eq. (16) that the circular frequency $\omega$ at the Hopf bifurcation points also depends exclusively on $\phi$ and $\tau_0$. This means that for the fixed values of $\phi$ and $\tau_0$ the frequency of oscillations does not depend on the value of the complex variable $c$. Since the real and imaginary parts of the complex variable $c_H$ at Hopf bifurcation point are expressed as:

$$c_{Hr} = \frac{\sqrt{1+\omega^2}}{2}\cos(\varphi) - \frac{1+\omega^2}{4}\cos(2\varphi) \qquad (17)$$

$$c_{Hi} = \frac{\sqrt{1+\omega^2}}{2}\sin(\varphi) - \frac{1+\omega^2}{4}\sin(2\varphi) \qquad (18)$$

the tangent of the argument of this constant $(k = c_{Hi}/c_{Hr})$ depends only on $\phi$. Hence, at each point $c(c_r; kc_r)$ the frequency of oscillations is the same, no matter the value of $c_r$. Let us assume, as an example, $\phi = 2.517$, which—according to Eqs. (16), (17) and (18)— corresponds to *ω* = 0.2292. Thus, the period of oscillations is equal to $T = \frac{2\pi}{\omega} = 27.417$ T=2πω=27.417. The numerical simulation of the time series |*z*(*τ*)| confirms the above analytical calculations, which is shown in Fig. 6. However, the result of simulation is only approximately true. Namely, at the increase of the simulation time one finds that this time series is unstable (Fig. 7). The consequence of it is the vanishing of the fractal structure of the image with the increase of simulation time. It should be added that the reverse phenomenon does not have to occur. This denotes that the removal of the time derivative from the complex differential equation with delay and thus the change of its form from the continuous form to the logistic one does not have to cause the generation of the fractal sets from this form. An example here is Eq. (8) analyzed in [5].

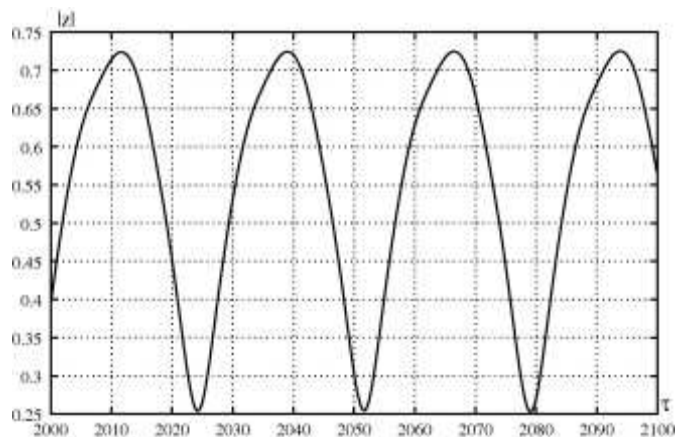

Fig. 6. Time story for small time simulation.

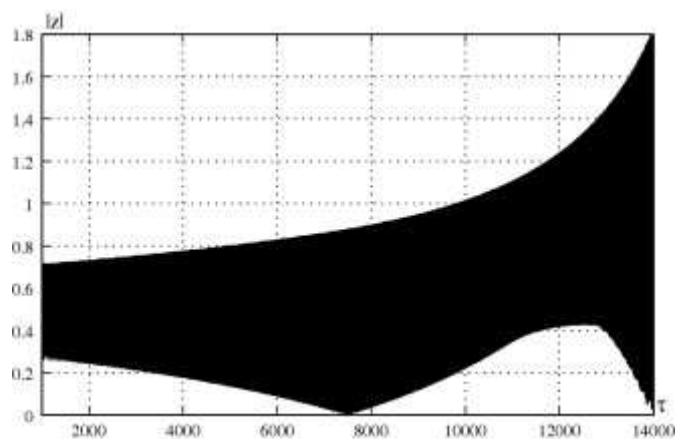

Fig. 7. Time story for large time simulation.

## 3. Concluding remarks

It has been shown that the logistic Mandelbrot's equation, extended by a dynamic term in the form of time derivative, does not generate the fractal images. This means that at real conditions, i.e. where any inertion occurs, the fractal structure generated by the extended Mandelbrot's equation is unstable. A following reflexion may be suggested. Commonly one observes that our environment consists of fractal structures. This results from the fact, that clouds, mountains, trees, snowflake, biological systems etc. are characterized by geometrical self-similarity. These structures have therefore definite fractional fractal dimensions. This feature concerns also the fractal generated by logistic Mandelbrot's equation (1) (Fig. 1) and also by extended Mandelbrot's equation (11) (Fig. 2). However, as was proved in the present work, the latter fractal vanishes in time. So, the question arises, does not the reality, which surrounds us, undergo the identical temporal transition? Does not tend it to an arranged form? Let us recall that on the atomic level the fractal structures do not exist.